\documentclass[12pt]{amsart}
\usepackage{fullpage}
\usepackage{amssymb,latexsym,amsmath,amscd,mathrsfs,amsthm}
\usepackage{graphics,graphicx}
\theoremstyle{plain}
\newtheorem{theorem}{Theorem}[section]
\newtheorem{corollary}[theorem]{Corollary}

\newcommand{\R}{\mathbb{R}}
\newcommand{\C}{\mathbb{C}}
\newcommand{\Q}{\mathbb{Q}}

\newcommand{\geom}[1]{\mathbb{#1}}

\newcommand{\D}{\mathrm{d}}
\newcommand{\id}{\mathrm{id}}
\newcommand{\gen}[1]{{\langle #1 \rangle}}
\newcommand{\conj}[1]{\overline{#1}}
\newcommand{\til}{\widetilde}

\renewcommand{\phi}{\varphi}
\renewcommand{\Im}{\mathrm{Im}\,}

\newcommand{\bdy}{\partial}
\newcommand{\Lie}[1]{\mathrm{#1}}

\newcommand{\Teich}{Teich\-m\"{u}l\-ler}
\newcommand{\dfn}[1]{\emph{#1}}

\begin{document}
\title[Orientation-reversing involutions of flat surfaces]%
{Orientation-reversing involutions of the genus 3 Arnoux--Yoccoz 
surface and related surfaces}\author{Joshua P. Bowman}
\address{Department of Mathematics, Cornell University, Ithaca, NY 14853}
\begin{abstract}
We present a new description of the genus 3 Arnoux--Yoccoz translation 
surface in terms of its Delaunay polygons and show that, up to affine 
equivalence, it belongs to two families of surfaces whose isometry groups 
include the dihedral group of the square.
\end{abstract}
\maketitle
\section{Introduction}\label{S:intro}
\subsection{Flat surfaces}\label{SS:fs}

%
%
%
%
%
%
%
%

For our purposes, a \dfn{flat surface} is a pair $(X,q)$, where $X$ is a 
Riemann surface and $q$ is a non-zero meromorphic quadratic differential 
of finite area on $X$. We will also speak of the flat surface $(X,\omega)$ 
instead of $(X,\omega^2)$ when $q = \omega^2$ is the global square of an 
abelian differential $\omega$; in this case, $(X,\omega)$ is also called a 
\dfn{translation surface}. A quadratic differential determines a canonical 
metric structure on the underlying surface (cf.\ \cite{jHhM79,cEfG97}); we 
will consider two flat surfaces to be the same when there is an isometry 
between them. 

Recall that $\Lie{SL}_2(\R)$ and $\Lie{GL}_2(\R)$ act on the space of 
translation surfaces, while $\Lie{PSL}_2(\R)$ and $\Lie{PGL}_2(\R)$ act 
on the space of flat surfaces. If $A \cdot (X,\omega) = (X,\omega)$ or 
$[A] \cdot (X,q) = (X,q)$, then $\det A = \pm 1$. The set of all 
$[A] \in \Lie{PGL}_2(\R)$ such that $[A]\cdot(X,q) = (X,q)$ form the 
\dfn{projective generalized Veech group} of $(X,q)$, which we will simply 
call the Veech group (cf.\ \cite{wV89,fHgS07}).

The results in this article began with a study of Delaunay polygons on the 
surface described in \S\ref{SS:ay}, and so we recall their definition 
(cf.\ \cite{hMjS91,iR94,wV97,aBbS07}).
A \dfn{Delaunay triangle} on $(X,q)$ is the image of a $2$-simplex on $X$, 
embedded on its interior, whose vertices lie in the set of cone points of 
$q$, whose edges are geodesic with respect to the metric of $q$, and whose 
(possibly immersed) circumcircle contains no cone points of $q$ in its 
interior.
A \dfn{Delaunay triangulation} of $(X,q)$ is a cell structure on $X$ whose 
$2$-cells are Delaunay triangles.
A \dfn{Delaunay polygon} is a Delaunay triangle or the union of two or 
more adjacent Delaunay triangles that share the same circumcircle.

A generic flat surface has a unique Delaunay triangulation. When it is not 
unique, we can start with any Delaunay triangulation and join two or more 
adjacent Delaunay triangles to form a Delaunay polygon. After forming all 
of the maximal Delaunay polygons from a Delaunay triangulation, we obtain 
the \dfn{Delaunay decomposition} of $(X,q)$, which is invariant under 
isometries of $(X,q)$.

\subsection{The Arnoux--Yoccoz surface}\label{SS:ay}

In \cite{pAjcY81}, P.~Arnoux and J.-C.~Yoccoz introduced a family 
of pseudo-Anosov diffeomorphisms, one for each genus at least $3$. 
Hubert--Lanneau \cite{pHeL05} showed that none of these surfaces has a 
Veech group containing any parabolic elements. The genus $3$ example, 
also appearing in \cite{pA88}, has received much attention recently. 
Hubert--Lanneau--M\"oller \cite{pHeLmM07} showed that the relevant 
abelian differential has a second, independent pseudo-Anosov element in 
its stabilizer, and using techniques introduced by C.~McMullen \cite{ctM07} 
in the study of genus $2$ orbits they showed that the $\Lie{SL}_2(\R)$-orbit 
of the genus $3$ example is dense in the largest possible region of the 
moduli space of abelian differentials. Here we give a new description of 
the genus $3$ surface in terms of its Delaunay polygons (of which there 
are only two kinds, up to isometry) and very simple gluing instructions.

\begin{figure}
\includegraphics[scale=1]{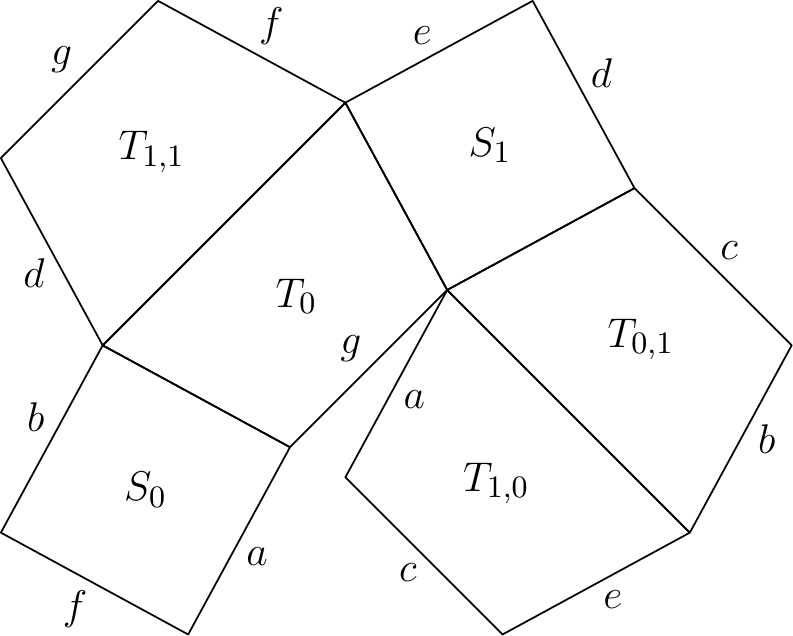}
\caption{The decomposition of 
$(X_\mathrm{AY},\omega_\mathrm{AY})$ into its Delaunay 
polygons. Edges with the same label are identified 
by translation.}
\label{F:AYquad}
\end{figure}

Let $\alpha \approx 0.543689$ be the real root of the polynomial 
$x^3 + x^2 + x - 1$. Let $S_0$ be the square with vertex set 
$\{ (0,0), (\alpha^2,\alpha), (\alpha^2-\alpha,\alpha^2+\alpha), 
(-\alpha,\alpha^2) \}$, and let $T_0$ be the trapezoid with vertex set 
$\{ (0,0), (1-\alpha,1-\alpha), (1-\alpha-\alpha^2,1), (-\alpha,\alpha^2) \}$. 
We form a flat surface $(X_\mathrm{AY},\omega_\mathrm{AY})$ from two 
copies of $S_0$ and four copies of $T_0$: reflecting $S_0$ across either 
a horizontal or vertical axis yields the same square $S_1$ (up to 
translation); we denote by $T_{1,0}$, $T_{0,1}$, and $T_{1,1}$ the 
reflections of $T_0$ across a vertical axis, across a horizontal axis, and 
across both, respectively. (In fact, $T_{1,0}$, $T_{0,1}$, and $T_{1,1}$ 
are all rotations of $T_0$ by multiples of $\pi/2$, but this description 
via reflections will be invariant under horizontal and vertical scaling, 
i.e., the \Teich\ geodesic flow.) Identify the long base of $T_0$ with 
the long base of $T_{1,1}$, as well as their short bases; do the same with 
$T_{1,0}$ and $T_{0,1}$. Each remaining side of a trapezoid is parallel to 
exactly one side of $S_0$ or $S_1$; identify by translations those sides 
which are parallel. (See Figure~\ref{F:AYquad}.)

The resulting flat surface $(X_\mathrm{AY},\omega_\mathrm{AY})$ has 
genus $3$ and two singularities each with cone angle $6\pi$. The images 
of $S_0$ and $T_0$ are the Delaunay polygons of $\omega_\mathrm{AY}$. 
$X_\mathrm{AY}$ is hyperelliptic; the hyperelliptic involution 
$\tau : X_\mathrm{AY} \to X_\mathrm{AY}$ is evident in Figure~\ref{F:AYquad} 
as rotation by $\pi$ around the centers of the squares and the midpoints 
of the edges joining two trapezoids; these six points together with the 
cone points are therefore the Weierstrass points of the surface. 
(See Figure~\ref{F:quotients} for the quotient of $X_\mathrm{AY}$ 
by $\tau$.) Moreover, $\omega_\mathrm{AY}$ is odd with respect to $\tau$, 
i.e., $\tau^*\omega_\mathrm{AY} = -\omega_\mathrm{AY}$.

\begin{figure}
\includegraphics[scale=1]{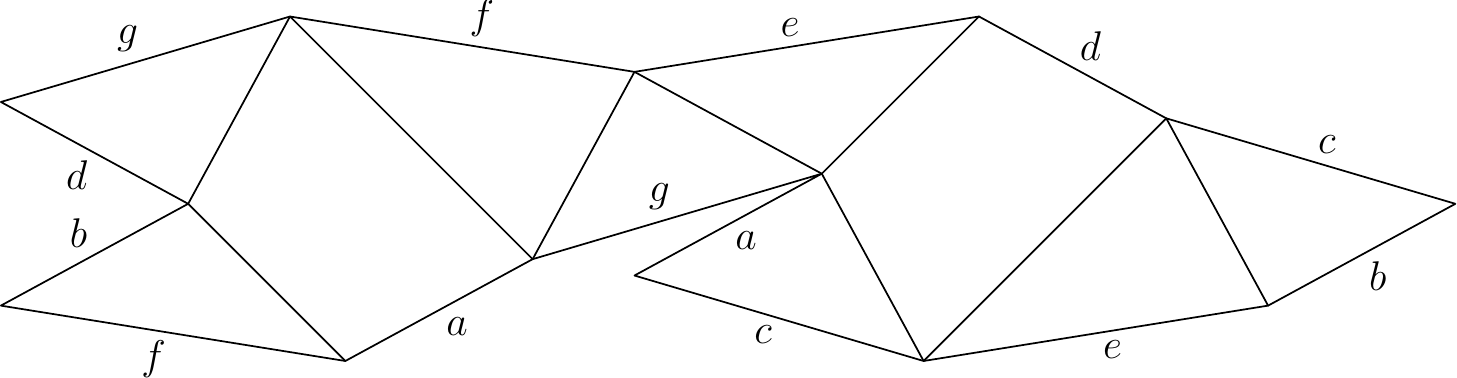}
\caption{The result of applying the Arnoux--Yoccoz 
pseudo-Anosov diffeomorphism to $\omega_\mathrm{AY}$. 
The original shapes of $S_0$ and $T_0$ and their copies 
can be reconstructed by matching edges.}\label{F:pA}
\end{figure}

The pseudo-Anosov diffeomorphism $\psi_\mathrm{AY}$ constructed by 
Arnoux--Yoccoz scales the surface by a factor of $1/\alpha$ in the 
horizontal direction and by $\alpha$ in the vertical direction. In 
Figure~\ref{F:pA} we show the result of applying this affine map to  
Figure~\ref{F:AYquad}, along with the new Delaunay edges. Two of the 
trapezoids---having the orientations of $T_{1,1}$ and $T_{1,0}$---are 
clearly visible; the squares and the other two trapezoids are constructed 
from the remaining triangles.

\subsection{$X_\mathrm{AY}$ as a cover of $\geom{RP}^2$}\label{SS:rp2}

The reflections applied to $S_0$ and $T_0$ in \S\ref{SS:ay} induce 
a pair of orientation-reversing involutions without fixed points on 
$X_\mathrm{AY}$. These can be visualized (as in Figure~\ref{F:AYquad}) 
as ``glide-reflections'', one along a horizontal axis and the other 
along a vertical axis. Both exchange $S_0$ and $S_1$. Let $\sigma_1$ be 
the involution that exchanges $T_0$ and $T_{1,0}$; i.e., its derivative 
is reflection across the horizontal axis. Let $\sigma_2$ be the 
involution that exchanges $T_0$ and $T_{0,1}$; i.e., its derivative is 
reflection across the vertical axis. The product of $\sigma_1$ and 
$\sigma_2$ is the hyperelliptic involution $\tau$, and neither sends any 
point of $X_\mathrm{AY}$ to its image by $\tau$. They therefore descend 
to a single involution $\sigma$ on $\geom{CP}^1$ without fixed points. 
The quotient of $\geom{CP}^1$ by $\sigma$ is homeomorphic to $\geom{RP}^2$.

\begin{figure}
\includegraphics[scale=1]{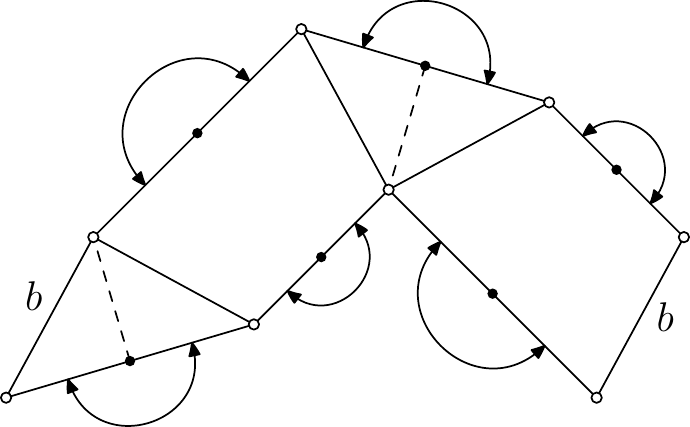}
\hspace{2cm}
\includegraphics[scale=1]{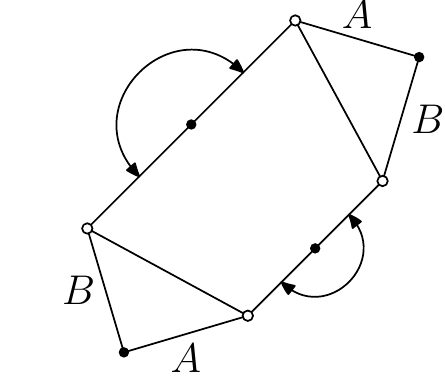}
\caption{Quotient surfaces of $X_\mathrm{AY}$.
{\sc left:} $\geom{CP}^1$ as the quotient of 
$X_\mathrm{AY}$ by $\gen{\tau}$. The edges marked $b$ 
are identified by translation.
{\sc right:} $\geom{RP}^2$ as the quotient of 
$X_\mathrm{AY}$ by the group $\gen{\sigma_1,\sigma_2}$. 
Edges with the same label are identified by a glide 
reflection along either a horizontal or a vertical axis. 
Dashed lines in the left picture indicate preimages of the 
segments labeled $B$ on the right.}\label{F:quotients}
\end{figure}

In fact, the presence of $\sigma_1$ and $\sigma_2$ is implicit in the 
work of Arnoux--Yoccoz. The original paper \cite{pAjcY81} begins with 
a measured foliation of $\geom{RP}^2$ with one ``tripod'' (a singular 
point of valence three) and three ``thorns'' (singular points of valence 
$1$), which is then lifted to the genus $3$ example we have described. 
In Figure~\ref{F:quotients}, right, we illustrate $\geom{RP}^2$ as the 
quotient of $X_\mathrm{AY}$ by the group generated by $\sigma_1$ and 
$\sigma_2$. In Figure~\ref{F:quotients}, left, we see $\geom{CP}^1$, on 
which $\sigma$ acts again by a ``glide-reflection'', which is the sheet 
exchange for the cover $\geom{CP}^1 \to \geom{RP}^2$. In both pictures 
we have drawn vertices that become tripods as open circles, and the 
vertices that become thorns as filled-in circles. The vertical foliation 
of the surface on the right of Figure~\ref{F:quotients} is the starting 
point of \cite{pAjcY81}.

In \S\ref{S:family} we will show that both 
$(X_\mathrm{AY},\omega_\mathrm{AY})$ and another affinely equivalent 
surface have \dfn{real structures} (orientation-reversing involutions 
whose fixed-point set is $1$-dimensional) that are not evident in the 
original construction. These additional structures will allow us to 
write equations for the surfaces and fit them into families of flat 
surfaces with a common group of isometries. In \S\ref{S:g2} we will 
transfer these results to genus $2$ quadratic differentials. In 
\S\ref{F:end} we will conclude by showing that we have found all the 
surfaces that are obtained by applying the geodesic flow to 
$(X_\mathrm{AY},\omega_\mathrm{AY})$ and have real structures.

\section{Two families of surfaces}\label{S:family}

\subsection{Labeling the Weierstrass points of $X_\mathrm{AY}$}\label{SS:wp}

As before, we denote the hyperelliptic involution of $X_\mathrm{AY}$ 
by $\tau$, and we let $\sigma_1$ and $\sigma_2$ be the involutions 
described in \S\ref{SS:rp2}, with $\sigma : \geom{CP}^1 \to \geom{CP}^1$ 
the involution covered by both $\sigma_1$ and $\sigma_2$. 

The purpose of this section is to show the following.

\begin{theorem}\label{T:family1}
The surface $(X_\mathrm{AY},\omega_\mathrm{AY})$ belongs to a family 
$(X_{t,u},\omega_{t,u})$, with $t > 1$ and $u > 0$, such that $X_{t,u}$ 
has the equation 
\begin{equation}\label{Eq:Xeqn}
y^2 = x (x - 1) (x - t) (x + u) (x + tu) (x^2 + tu),
\end{equation}
and $\omega_{t,u}$ is a multiple of $x\,\D{x}/y$ on $X_{t,u}$.
\end{theorem}

Each of the surfaces in Theorem~\ref{T:family1} has a pair of real 
structures $\rho_1$ and $\rho_2$ whose product is again $\tau$, 
and which therefore descend to a single real structure $\rho$ on 
$\geom{CP}^1$. Any product of the form $\rho_i \sigma_j$ ($i,j \in \{1,2\}$) 
is a square root of $\tau$, and therefore the group generated by 
$\{\sigma_1,\sigma_2,\rho_1,\rho_2\}$ is the dihedral group of the square. 
We will exhibit these isometries in our presentation of 
$(X_\mathrm{AY},\omega_\mathrm{AY}$). In \S\ref{SS:add} we will look 
at surfaces in this family that have additional symmetries.

Let $\varpi : X_\mathrm{AY} \to \geom{CP}^1$ be the degree 2 map induced 
by $\tau$, i.e., $\varpi\circ\tau = \varpi$. We can normalize $\varpi$ 
so that the zeroes of $\omega_\mathrm{AY}$ are sent to $0$ and $\infty$, 
and the midpoint of the short edge between $T_0$ and $T_{1,1}$ is sent 
to $1$. We wish to find the images of the remaining Weierstrass points, 
so that we can write an affine equation for $X_\mathrm{AY}$ in the form 
$y^2 = P(x)$, where $P$ is a degree 7 polynomial with roots at $0$ and 
$1$. Hereafter we assume that $\varpi$ is the restriction to 
$X_\mathrm{AY}$ of the coordinate projection $(x,y) \mapsto x$. 
Consequently, we may consider each Weierstrass point as either a point 
$(w,0)$ that solves $y^2 = P(x)$ or simply as a point $w$ on the $x$-axis.

Each of the real structures $\rho_1$ and $\rho_2$ has a fixed-point set 
with three components: in one case, say $\rho_1$, the real components are 
the line of symmetry shared by $T_0$ and $T_{1,1}$, and the two bases of 
$T_{1,0}$ and $T_{0,1}$. The fixed-point set of $\rho_2$ is then the union 
of the corresponding lines in the orthogonal direction. Because $\rho_1$ 
and $\rho_2$ fix the points $0$, $1$, and $\infty$, $\rho$ fixes the real 
axis; therefore $\rho$ is simply complex conjugation.

With this normalization, the involution $\sigma$ on $\geom{CP}^1$ 
exchanges $0$ and $\infty$ and preserves the real axis; therefore 
$\sigma$ has the form $x \mapsto -r/\conj{x}$ for some real $r > 0$.

Let $s = (s,0)$ be the center of $S_0$. Then 
$\rho_1(s) = \rho_2(s) = \sigma_1(s) = \sigma_2(s)$ is the center of $S_1$, 
which implies $\rho(s) = \sigma(s)$, i.e., $\conj{s} = -r/\conj{s}$. The 
solutions to this equation are $\pm i\sqrt{r}$. By considering the location 
of the fixed-point sets of $\rho_1$ and $\rho_2$, we see that the image of 
$S_0$ by $\varpi$ lies in the upper half-plane; therefore $s = i\sqrt{r}$, 
and $-i\sqrt{r}$ is the center of $S_1$.

Let $t$ be the midpoint of the long edge of $T_0$. Applying $\sigma_1$ or 
$\sigma_2$ shows that the midpoint of the long edge of $T_{1,0}$ is at $-r/t$. 

We already know that $1$ is the center of the short edge of $T_0$. Since 
the short edge of $T_{0,1}$ is the image of this edge by $\sigma_1$ or 
$\sigma_2$, the midpoint of this edge must be at $-r$.

\begin{figure}
\includegraphics[scale=1]{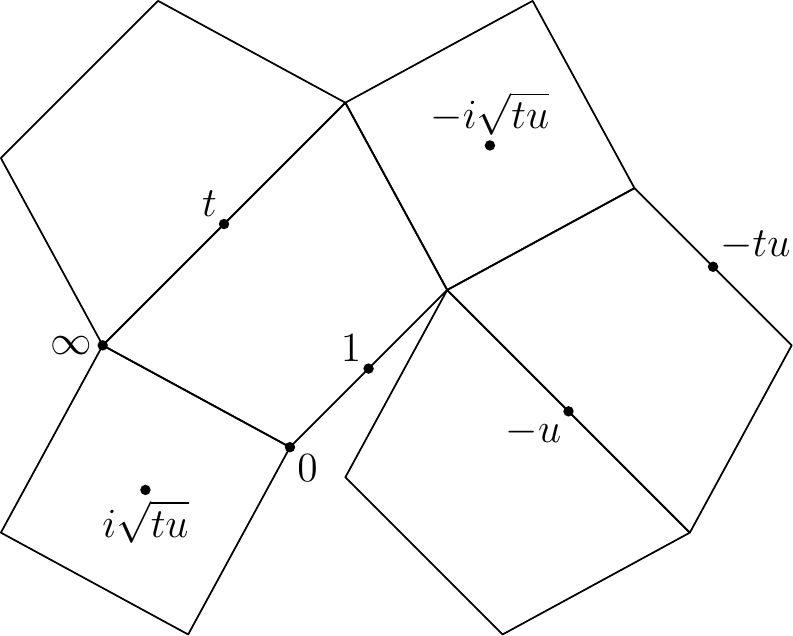}
\caption{The Weierstrass points of $X_\mathrm{AY}$, 
following normalization ($t > 1$, $u > 0$). The real 
structures $\rho_1$ and $\rho_2$ appear as reflections 
in the lines of slope $\pm1$.}\label{F:AYWP}
\end{figure}

To simplify notation, let us make the substitution $u = r/t$, so that 
$r = tu$ (hence $\sigma$ has the form $\sigma(x) = -tu/\conj{x}$). Thus 
$X_\mathrm{AY}$ has the equation~\eqref{Eq:Xeqn} for some $(t,u) = 
(t_\mathrm{AY},u_\mathrm{AY})$. Furthermore, $\omega_\mathrm{AY}$ is the 
square root of a quadratic differential on $\geom{CP}^1$ with simple zeroes 
at $0$ and $\infty$ and simple poles at $1$, $t$, $-u$, $-tu$, and 
$\pm i\sqrt{tu}$. There is therefore some complex constant $c$ such that 
\[
{\omega_\mathrm{AY}}^2 
= \varpi^* 
  \left( 
  \frac{cx}{(x - 1) (x - t) (x + u) (x + tu) (x^2 + tu)}\, \D{x}^2 
  \right),
\]
i.e., $\omega_\mathrm{AY} = \pm\sqrt{c}\,x\,\D{x}/y$. This 
establishes Theorem~\ref{T:family1}.

\subsection{Integral equations}\label{SS:int}

To find $t_\mathrm{AY}$ and $u_\mathrm{AY}$ requires solving a system 
of equations involving hyperelliptic integrals, which we establish in 
this section using relative periods of $\omega_\mathrm{AY}$. Choose a 
square root of 
\[
f_{t,u}(x) = \frac{x}{(x-1)(x-t)(x+u)(x+tu)(x^2+tu)}
\]
in the open first quadrant such that its extension $\sqrt{f_{t,u}(x)}$ 
to the complement of $\{1,t,i\sqrt{tu}\}$ in the closed first quadrant 
is positive on the open interval $(0,1)$. Let $\eta_0$ be the Delaunay 
edge between $S_0$ and $T_0$; $\varpi(\eta_0)$ is then a curve from $0$ 
to $\infty$ in the first quadrant. Integrating $\sqrt{cf_{t,u}(x)}\,\D{x}$ 
on the portion of the first quadrant below $\varpi(\eta_0)$ will then give 
a conformal map to half of $T_0$. We will be interested in integrals along 
the real axis.

The vector from $0$ to $1$ along the short side of $T_0$ is 
$\frac{1}{2}(1-\alpha)(1 + i)$, while the the line of symmetry of 
$T_0$ from $1$ to $t$ gives the vector $\frac{1}{2}(\alpha+\alpha^2)(-1+i)$. 
Observe that 
\[
i \cdot (1 - \alpha) (1 + i) = \alpha \cdot (\alpha + \alpha^2) (-1 + i),
\]
and therefore
\begin{equation}\label{Eq:int1}
i \int_0^1 \sqrt{cf_{t,u}(x)}\,\D{x} 
= \alpha \int_1^t \sqrt{cf_{t,u}(x)}\,\D{x}.
\end{equation}
Similarly, the vector from $t$ to $\infty$ along the long side of 
$T_0$ is $\dfrac{1}{2} (1 - \alpha^2) (-1 - i)$, and because 
$1 - \alpha^2 = (1 + \alpha) (1 - \alpha)$, we have 
\begin{equation}\label{Eq:int2}
-(1 + \alpha) \int_0^1 \sqrt{cf_{t,u}(x)}\,\D{x} 
= \int_t^\infty \sqrt{cf_{t,u}(x)}\,\D{x}.
\end{equation}
In both equations we can cancel out the $c$, which was ever only 
a global (complex) scaling factor anyway. Now bring $i$ under the 
square root on the right-hand side of \eqref{Eq:int1} in order to 
make the radicand positive. We thus obtain from \eqref{Eq:int1} and 
\eqref{Eq:int2} the system of (real) integral equations
\begin{equation}\label{Eq:hypellsystem}
\begin{cases}
\displaystyle
\phantom{(1+\alpha)}\int_0^1 \sqrt{f_{t,u}(x)}\,\D{x} 
= \alpha \int_1^t \sqrt{-f_{t,u}(x)}\,\D{x} \\
\displaystyle
(1+\alpha) \int_0^1 \sqrt{f_{t,u}(x)}\,\D{x} 
= -\int_t^\infty \sqrt{f_{t,u}(x)}\,\D{x}
\end{cases}
\end{equation}
whose solution is the desired pair $(t_\mathrm{AY},u_\mathrm{AY})$. 
Using numeric methods, we find 
\[
t_\mathrm{AY} \approx 1.91709843377 \qquad\mbox{and}\qquad
u_\mathrm{AY} \approx 2.07067976690.
\]
We conjecture that $t_\mathrm{AY}$ and $u_\mathrm{AY}$ lie in some 
field of small degree over $\Q(\alpha)$.

\subsection{Other exceptional surfaces in this family}\label{SS:add}

An examination of the geometric arguments in \S\ref{SS:wp} and an 
application of the principle of continuity to $t$ and $u$ show the 
following: 

\begin{theorem}\label{T:trapezoid}
Every $(X_{t,u},\omega_{t,u})$ as in Theorem~\ref{T:family1} can be 
formed by replacing $T_0$ in the description from \S\ref{SS:ay} with 
an isosceles trapezoid $T$, $S_0$ with the square built on a leg of 
$T$, and the copies of $T_0$ with the rotations of $T$ by $\pi/2$.
\end{theorem}

The placement of $t$ and $u$ on $\R$ determines the shape of the 
trapezoid $T$, and any isosceles trapezoid may be obtained by an 
appropriate choice of $t$ and $u$. In this section, we examine 
certain shapes that give $(X_{t,u},\omega_{t,u})$ extra symmetries 
and determine the corresponding values of $t$ and $u$. We continue 
to use $\tau$ to denote the hyperelliptic involution of $X_{t,u}$.

Suppose that $T$ is a rectangle. Then there are two orthogonal closed 
trajectories, running parallel to the sides of $T$ and connecting the 
centers $\pm i\sqrt{tu}$ of the squares, and either of these can be 
made into the fixed-point set of a real structure on $X_{t,u}$. The 
product of these two real structures is again $\tau$, so they descend 
to a single real structure on $\geom{CP}^1$. This real structure 
exchanges $0$ with $\infty$ and fixes $\pm i\sqrt{tu}$, so it must be 
inversion in the circle $|x|^2 = tu$. It also exchanges $1$ with $t$, 
which implies $1\cdot t = tu$, i.e., $u = 1$. The remaining parameter 
$t$ is determined by solving the single integral equation 
\[
\int_0^1 \sqrt{\frac{x}{(x^2 - 1) (x^2 - t^2) (x^2 + t)}} \,\D{x} = 
\mu \int_1^t \sqrt{\frac{-x}{(x^2 - 1) (x^2 - t^2) (x^2 + t)}} \,\D{x}
\]
where $2\mu$ is the ratio of the width of $T$ to its height. Recall 
that an \dfn{origami}, also called a \dfn{square-tiled surface}, is a 
flat surface that covers the square torus with at most one branch point 
(cf.\ \cite{gS00,aEaO01,aZ02}). By looking at rational values of $\mu$, 
we have the following result:

\begin{corollary}
The family $(X_{t,1},\omega_{t,1})$ contains a dense set of origamis.
\end{corollary}

These are not the only $(X_{t,u},\omega_{t,u})$ that are origamis, 
however. If $T$ is a trapezoid whose legs are orthogonal to each other, 
then $(X_{t,u},\omega_{t,u})$ is again an origami.

\subsection{Second family of surfaces}

\begin{figure}
\includegraphics[scale=1]{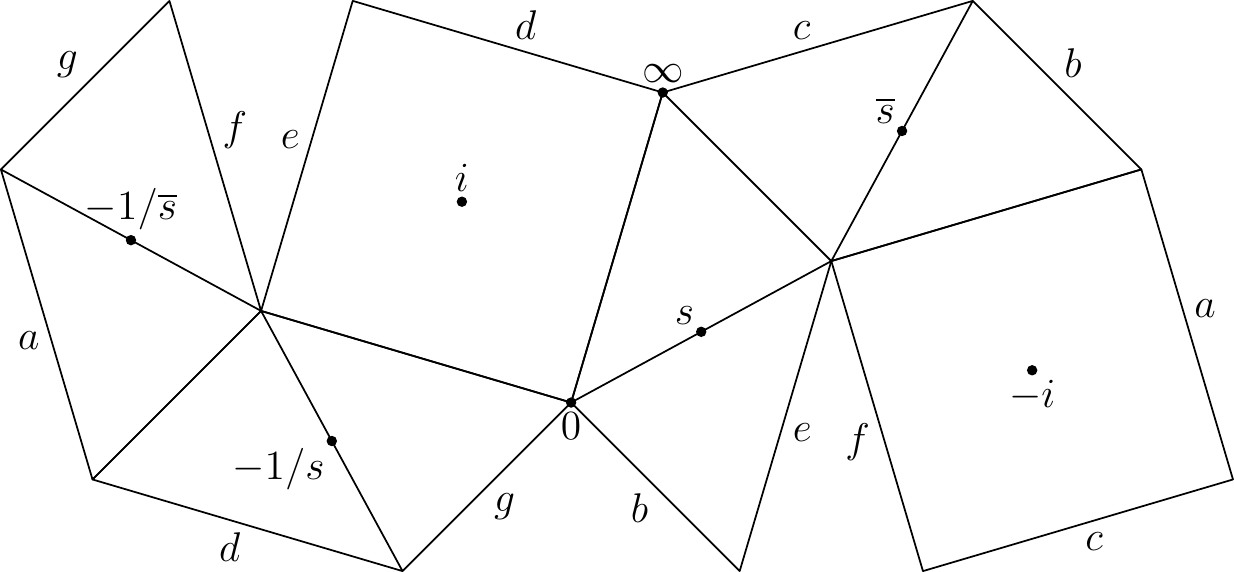}
\caption{Another surface in the $\Lie{GL}_2(\R)$-orbit 
of $\omega_\mathrm{AY}$ with additional real structures. 
Edges with the same label are identified.}\label{F:AY2}
\end{figure}

Conjugating $\rho_1$ by the pseudo-Anosov element $\psi_\mathrm{AY}$ 
guarantees the existence of another orientation-reversing involution in 
the affine group of $\omega_\mathrm{AY}$. This element fixes a point 
``half-way'' (in the Teichm\"uller metric, for instance) between 
$\omega_\mathrm{AY}$ and its image by $\psi_\mathrm{AY}$, lying in the 
Teichm\"uller disk of $(X_\mathrm{AY},\omega_\mathrm{AY})$. This surface 
can be found either by scaling the vertical direction by $\sqrt{\alpha}$ 
and the horizontal direction by $1/\sqrt{\alpha}$ or, to keep our 
coordinates in the field $\Q(\alpha)$, just by scaling the horizontal by 
$1/\alpha$. This surface, which we will denote 
$(X_\mathrm{AY}',\omega_\mathrm{AY}')$, is shown in Figure~\ref{F:AY2}, 
along with its Delaunay polygons.

\begin{theorem}\label{T:family2}
The surface $(X'_\mathrm{AY},\omega'_\mathrm{AY})$ belongs to a 
family $(X_s,\omega_s)$, with $\Im s > 0$ and $s \ne i$, such that 
$X_s$ has the equation 
\begin{equation}\label{Eq:X2eqn}
y^2 = x (x^2 + 1) (x - s) (x - \conj{s}) (x + 1/s) (x + 1/\conj{s}),
\end{equation}
and $\omega_s$ is a multiple of $x\,\D{x}/y$ on $X_s$.
\end{theorem}

Again, we have two real structures $\rho'_1$ and $\rho'_2$ whose 
product is the hyperelliptic involution $\tau$. Each of these only 
has one real component, however: the union of the sides of the 
parallelograms running parallel to the axis of reflection. The only 
Weierstrass points that lie on these components are $0$ and $\infty$; 
the remaining Weierstrass points are the centers of the squares and 
of the parallelograms. We again let $\rho'$ be the induced real 
structure on $\geom{CP}^1$ and assume that it fixes the real axis 
(this we can do because we have only fixed the positions of two 
points on $\geom{P}^1$), so that the remaining Weierstrass points 
come in conjugate pairs.

The fixed-point free involutions $\sigma_1$ and $\sigma_2$ from 
\S\ref{SS:rp2} again preserve the union of the real loci of $\rho'_1$ 
and $\rho'_2$, and therefore they descend to a fixed-point free 
involution $\sigma$ of the form $x \mapsto -r/\conj{x}$, with $r$ real 
and positive. We have one more free real parameter for normalization, 
so we can assume $r = 1$. This implies that the centers of the squares 
are at $\pm i$. Let $s$ be the center of one of the parallelograms; 
then applying $\rho'_1$ and $\sigma_1$ shows that the remaining 
Weierstrass points are $\conj{s}$, $1/s$, and $1/\conj{s}$. Using 
developing vectors again, we can find equations that define $s$, in a 
manner analogous to finding \eqref{Eq:hypellsystem}. 

As an analogue to Theorem~\ref{T:trapezoid}, we have:

\begin{theorem}
Every $(X_s,\omega_s)$ as in Theorem~\ref{T:family2} can be formed 
from a parallelogram $P$, a square built on one side of $P$, the 
rotation of $P$ by $\pi/2$, and the images of $P$ and its rotation 
by reflection across their remaining sides.
\end{theorem}

The shape of $P$ is determined by the value of $s$. If $s = 
\frac{1}{2}(\sqrt{3}+i)$, then $P$ becomes a square, and we obtain 
one of the ``escalator'' surfaces in \cite{sLrS07}. More generally, 
if $s$ is any point of the unit circle, then $P$ is a rectangle, 
and inversion in the unit circle corresponds to another pair of 
real structures on $X$, which are the reflections across the axes 
of symmetry of $P$. By considering those rectangular $P$ whose side 
lengths are rationally related, we have as before:

\begin{corollary}
The family $(X_{e^{i\theta}},\omega_{e^{i\theta}})$ (with 
$0 < \theta < \pi/2$) contains a dense set of origamis.
\end{corollary}

Another origami appears when $P$ is composed of a pair of right 
isosceles triangles so that $s$ lies not on the hypotenuse, but on 
a leg of each.

\section{Quadratic differentials and periods on genus $2$ surfaces}\label{S:g2}

We do not know how to compute the rest of the periods for $X_{t,u}$ or 
$X_s$, apart from those of $\omega_{t,u}$ or $\omega_s$, respectively. 
In this section, however, we consider the periods of certain related 
genus $2$ surfaces, which demonstrate remarkable relations.

Let $X$ be any hyperelliptic genus $3$ surface with an abelian differential 
$\omega$ that is odd with respect to the hyperelliptic involution and has 
two double zeroes. The pair $(X,\omega)$ has a corresponding pair $(\Xi,q)$, 
where $\Xi$ is a genus $2$ surface and $q$ is a quadratic differential on 
$\Xi$ with four simple zeroes. Geometrically, the correspondence may be 
described as follows: two of the zeroes of $\omega$ are at Weierstrass 
points of $X$, hence $(X,\omega^2)$ covers a flat surface 
$(\geom{CP}^1,\til{q})$ where $\til{q}$ has six poles and two simple zeroes 
(Figure~\ref{F:quotients}). Then $(\Xi,q)$ is the double cover of 
$(\geom{CP}^1,\til{q})$ branched at the poles of $\til{q}$. In our cases, 
the genus $2$ surface may be obtained by cutting along opposite sides of 
one of the squares in Figure \ref{F:AYquad} or \ref{F:AY2}, then regluing 
each of these via a rotation by $\pi$ to the free edge provided by cutting 
along the other (cf.\ \cite{eL04,sV05}).

First we consider the family $(X_{t,u},\omega_{t,u})$ and the related 
genus $2$ flat surfaces $(\Xi_{t,u},q_{t,u})$. To be explicit, the 
defining expressions of both types of surfaces are:
\begin{align*}
X_{t,u} :\ 
& y^2 = x (x - 1) (x - t) (x + u) (x + tu) (x^2 + tu), 
&& \omega_{t,u} = \frac{x\,\D{x}}{y}; \\
\Xi_{t,u} :\
& y^2 = (x - 1) (x - t) (x + u) (x + tu) (x^2 + tu),
&& q_{t,u} = \frac{x\,\D{x}^2}{y^2}.
\end{align*}
The order $4$ rotation $\rho_1 \sigma_1$ of $X_{t,u}$ persists on 
$\Xi_{t,u}$. Following R.~Silhol \cite{rS06}, we find a new parameter $a$, 
depending on $t$ and $u$, so that the Riemann surface 
\[
\Xi_a : y^2 = x (x^2 - 1) (x - a) (x - 1/a)
\]
is isomorphic to $\Xi_{t,u}$. Doing so simply requires a change of 
coordinates in $x$, namely 
\[
\Phi(x) = i \sqrt{tu} \frac{(x - 1)}{(x + tu)}.
\]
Then $\Phi(1) = 0$, $\Phi(-tu) = \infty$, and $\Phi(\pm i\sqrt{tu}) = \mp 1$. 
The images of $t$ and $u$ by $\Phi$ are 
\[
a = i \sqrt{\frac{u}{t}} \frac{(t - 1)}{(u + 1)} 
\qquad\mbox{and}\qquad
\frac{1}{a} = i \sqrt{\frac{t}{u}} \frac{(1 + u)}{(1 - t)}. 
\]
Because $t > 1$ and $u > 0$, $a$ lies on the positive imaginary axis 
and $1/a$ lies on the negative imaginary axis. The involution $\rho$ 
becomes reflection across the imaginary axis. The images of $0$ and 
$\infty$ by $\Phi$ are 
\[
\Phi(0) = \frac{i}{\sqrt{tu}} 
\qquad\mbox{and}\qquad
\Phi(\infty) = \frac{\sqrt{tu}}{i},
\]
so the image of $q_{t,u}$ on $\Xi_a$ is a scalar multiple of 
\[
\left( x - \frac{i}{\sqrt{tu}} \right) \left( x - \frac{\sqrt{tu}}{i} \right) 
	\frac{\D{x}^2}{y^2}
= \left( x^2 + i\left(\frac{tu - 1}{\sqrt{tu}}\right)x + 1 \right)
	\frac{\D{x}^2}{y^2}
\]
These calculations imply that, for each pair $(t_0,u_0)$, there is a 
one-parameter family of surfaces $(\Xi_{t,u},q_{t,u})$ such that 
$\Xi_{t,u}$ is isomorphic to $\Xi_{t_0,u_0}$ while $q_{t,u}$ and 
$q_{t_0,u_0}$ represent different differentials on the abstract Riemann 
surface.

Now we apply the same analysis to the second family. This time we are 
moving from $(X_s,\omega_s)$ to $(\Sigma_s,q_s)$, as defined below:
\begin{align*}
X_s :\ 
& y^2 = x (x^2 + 1) (x - s) (x - \conj{s}) (x + 1/s) (x + 1/\conj{s}), 
&& \omega_s = \frac{x\,\D{x}}{y}; \\
\Sigma_s :\
& y^2 = (x^2 + 1) (x - s) (x - \conj{s}) (x + 1/s) (x + 1/\conj{s}),
&& q_s = \frac{x\,\D{x}^2}{y^2}.
\end{align*}
We change coordinates in $x$ using 
\[
\Psi(x) = i\left(\frac{x - s}{sx + 1}\right)
\]
so that $\Psi(s) = 0$, $\Psi(-1/s) = \infty$, and $\Psi(\pm i) = \mp 1$. 
This time we get the curve $y^2 = x (x^2 - 1) (x - a) (x - 1/a)$, where 
\[
a = \Phi(\conj{s}) = \frac{2\,\Im s}{1 + |s|^2}
\qquad\mbox{and}\qquad
\frac{1}{a} 
= \Phi\left(-\frac{1}{\conj{s}}\right) 
= \frac{1 + |s|^2}{2\,\Im s}.
\]
Here we have $0 < a < 1$ and $1/a > 1$; $\rho'$ becomes inversion 
in the unit circle. The points $0$ and $\infty$ on $\Sigma_s$ become 
$\Phi(0) = -is$ and $\Phi(\infty) = i/s$. Again, we find just a 
one-parameter family of genus $2$ Riemann surfaces, each carrying a 
one-parameter family of quadratic differentials corresponding to 
distinct surfaces $X_s$.

In \cite{rS06}, it is shown that the full period matrix for any of 
the surfaces $\Xi_a$ can be expressed in terms of a single parameter, 
thanks to the fourfold symmetry of the surface. This parameter is the 
ratio of $\int_{-1}^0 \phi$ and $\int_0^{1/a} \phi$, where
$\phi = \frac{\D{x}}{y} - \frac{x\,\D{x}}{y}$. This ratio is real 
precisely when $a$ lies on the positive imaginary axis, as in our 
first family, and in these cases the period matrix of $\Xi_a$ is 
purely imaginary.

\section{Final remarks}\label{F:end}

\begin{figure}
\includegraphics{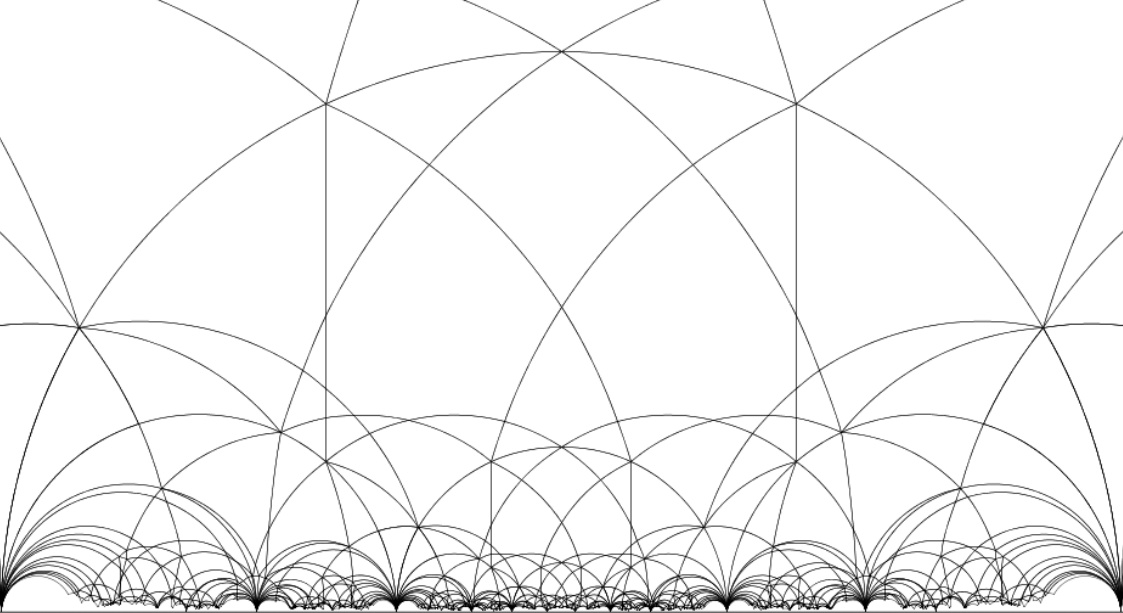}
\caption{The iso-Delaunay tessellation of $\geom{H}$ arising from 
$(X_\mathrm{AY},\omega_\mathrm{AY})$}\label{F:AYidt}
\end{figure}

The involutions we have exhibited also act on the \Teich\ disk generated 
by $(X_\mathrm{AY},\omega_\mathrm{AY})$, and their effects can be seen 
via the \dfn{iso-Delaunay tessellation} shown in Figure~\ref{F:AYidt}. 
The open regions in this picture correspond to combinatorial classes of 
Delaunay triangulations of surfaces in the $\Lie{SL}_2(\R)$-orbit of 
$(X_\mathrm{AY},\omega_\mathrm{AY})$; points of the $1$-skeleton 
correspond to surfaces with two or more Delaunay triangulations. Because 
Delaunay triangulations are not changed by any rotation of the surface, 
this picture can be drawn in the upper half-plane $\geom{H}$ rather than 
its unit tangent bundle. Iso-Delaunay tessellations have been described 
previously in \cite{jpb06} and \cite{wV08}. We will not define them here, 
but simply illustrate how elements of the generalized Veech group 
$\Gamma \subset \Lie{GL}_2(\R)$ of $(X_\mathrm{AY},\omega_\mathrm{AY})$ 
may be seen to act on the tessellation in Figure~\ref{F:AYidt}. (The 
hyperelliptic involution, having derivative $-\id$, acts trivally.)

Each element of $\Gamma$ acts on $\geom{H}$ by an isometry, preserving 
or reversing orientation according to the sign of its determinant. 
Figure~\ref{F:AYidt} is symmetric with respect to the central axis 
(the imaginary axis in $\C$); both $\sigma_1$ and $\sigma_2$ yield 
elements of $\Gamma$ that reflect $\geom{H}$ across this axis. The 
hyperbolic element of $\Gamma$ corresponding to $\psi_\mathrm{AY}$ 
fixes the points $0$ and $\infty$ in $\bdy\geom{H}$ and translates 
points along the imaginary axis by $z \mapsto z/\alpha^2$. A sequence 
of concentric circles is visible in the tessellation; these are 
related by $\psi_\mathrm{AY}$, and one is the unit circle, so their 
radii are all powers of $1/\alpha^2 \approx 3.38$.

There are two kinds of distinguished points on the central axis: 
ones where two geodesics meet and ones where three geodesics meet. 
The latter are those whose corresponding surface is isometric to 
$(X_\mathrm{AY},\omega_\mathrm{AY})$, while the former correspond 
to $(X'_\mathrm{AY},\omega'_\mathrm{AY})$. The real structures 
$\rho_1$ and $\rho_2$ (resp.\ $\rho'_1$ and $\rho'_2$) yield an 
element of $\Gamma$ that reflects $\geom{H}$ across the unit circle 
(resp.\ across the circle $|z| = 1/\alpha$). The order $4$ rotations 
of $(X_\mathrm{AY},\omega_\mathrm{AY})$ and 
$(X'_\mathrm{AY},\omega'_\mathrm{AY})$ are thus visible as the order 
$2$ rotations of $\geom{H}$ around these distinguished points.

If any other flat surface on the central axis had real structures, 
then its symmetries, too, would have to induce a reflection of 
$\geom{H}$ that preserves the tessellation. No such point exists; 
therefore we have described all the surfaces within the orbit of 
$(X_\mathrm{AY},\omega_\mathrm{AY})$ under the geodesic flow that 
demonstrate additional symmetries.

\bibliographystyle{math}
\bibliography{ABC}

\end{document}